\documentclass[10pt]{article}

\usepackage[top=1in,bottom=1in,left=1.25in,right=1.25in]{geometry}

\usepackage[absolute,overlay]{textpos}

\usepackage{graphicx}
\usepackage{ragged2e}
\usepackage{tabto}

\usepackage{amsmath}
\usepackage{amsthm}
\usepackage{amsfonts}
\usepackage{amssymb}

\usepackage{float}
\usepackage{multirow}
\usepackage{tikz}
\usepackage{graphics}
\usepackage{color}
\usepackage{url}
\usepackage{afterpage} 
\usepackage{multicol}
\usepackage[font=footnotesize]{caption}
\usepackage[font=footnotesize]{subcaption}
\usepackage{marginnote}
\usepackage{tcolorbox}

\usepackage{lmodern}

\usepackage[sf,bf,medium]{titlesec}
\usepackage{pgfplots}

\pgfplotsset{compat=1.18}

\usepackage{bm}
\usepackage{booktabs}
\usepackage{siunitx}
\usepackage{isomath}

\usepackage[plainpages=false, colorlinks=true, citecolor=blue, filecolor=blue,
   linkcolor=blue, urlcolor=blue]{hyperref}
\usepackage{cleveref}

\usepackage{multicol}
\newcommand{\opnm}{\operatorname}
\newcommand{\opd}{{\rm d}}

\newcommand{\lp}{\left(}
\newcommand{\rp}{\right)}

\usepackage[sort,compress]{cite}

\newcommand\bbR{\mathbb R}

\newcommand\bx{\boldsymbol{x}}

\newcommand\by{\boldsymbol{y}}

\newcommand\bn{\boldsymbol n}
\newcommand\btau{\boldsymbol \tau}

\newcommand\dd{{\rm{d}}}

\newtheorem{remark}{\sffamily Remark}

\newtheorem{conjecture}{\sffamily Conjecture}

\newcommand{\cK}{\mathcal K}

\numberwithin{equation}{section}

\newcommand{\bs}{\boldsymbol}

\newcommand{\figref}[1]{\figurename~\hyperref[#1]{\ref{#1}}}
\newcommand{\tableref}[1]{\tablename~\hyperref[#1]{\ref{#1}}}

\renewcommand{\phi}{\varphi}

\binoppenalty=\maxdimen
\relpenalty=\maxdimen

\newcommand{\gammaunit}{{\gamma_0}}

\usepackage{setspace}
\usetikzlibrary{calc}
\usepackage{abstract}

\begin{document}

{\centering {\Large Integral equations for flexural scattering problems with periodic boundaries\\}
}

\begin{center}
{Fruzsina Agocs\footnote{Department of Computer Science, University of Colorado Boulder.}, Tristan Goodwill\footnote{Department of Statistics and the Committee on Computational and Applied Mathematics, University of Chicago.}, Jeremy G. Hoskins\footnote{Department of Statistics and the Committee on Computational Applied Mathematics, University of Chicago, USA and NSF-Simons National Institute for Theory and Mathematics in Biology, Chicago, IL.}, and  Peter Nekrasov\footnote{Committee on Computational and Applied Mathematics, University of Chicago.}}
\end{center}

\begin{abstract}
We develop a method for computing the scattering of flexural waves off of a periodic wall or a periodic line of scatterers. These waves model the fluctuations of thin plates with periodic clamped, supported, or free edges. We use the Floquet--Bloch transform to convert the problem into a collection of uncoupled quasi-periodic problems. We then solve each quasi-periodic problem efficiently and accurately using a novel integral equation based on the quasi-periodic flexural Green's function. Finally, we show how the proposed method can be used to simulate scattering from junctions of semi-infinite lines of scatterers.

\end{abstract}
  

\section{Introduction}
\label{sec:intro}

In this paper, we consider the propagation of time-harmonic \emph{flexural waves} in thin plates with infinite periodic boundaries, including linear arrays of holes and clamped regions. Of particular interest are the existence of \emph{edge modes} which are localized near the boundary and propagate along it without decay. 
These modes were first demonstrated for free plates by Konenkov in 1960 \cite{konenkov1960rayleigh} and later, independently, by Sinha \cite{sinha1974some} and Thurston and McKenna \cite{thurstonmckenna}. This theory was also extended to free plates with periodic clamped pinning points along the boundary \cite{evans2008flexural}. More recently, flexural edge modes have appeared in a number of applications, \emph{inter alia} topological waves in elastic metamaterials \cite{wang2018topological,mousavi2015topologically,miniaci2018experimental}, scattering from an array of pins\cite{haslinger2016dynamic}, guided modes in plates made from LEGOs \cite{shen2024guided}, as well as linear arrays of quasiperiodically arranged resonators \cite{marti2021edge}.

As our model, we take the time-harmonic flexural wave equation
\begin{equation}\label{eq:flex}
    \Delta^2 u - k^4 u= 0, 
\end{equation}
in an unbounded region $\Omega \subset \mathbb{R}^2$ with a periodic boundary $\gamma$. On $\gamma$, we impose either \emph{clamped plate} boundary conditions
\begin{equation}
        u|_\gamma = f \quad\text{and}\quad \partial_{\bn} u|_\gamma = g ,
\end{equation}
where $\bn$ is the unit normal to $\gamma$ pointing into $\Omega$, 
\emph{supported plate} boundary conditions
\begin{equation}
    u|_\gamma = f \quad\text{and}\quad\nu \Delta u +(1-\nu) \partial_{\bn}^2u|_\gamma = f ,
\end{equation}
where $\nu\in[-1,1/2]$ is Poisson's ratio \cite{landau2012theory}, or \emph{free plate} boundary conditions
\begin{equation}
    \nu \Delta u +(1-\nu) \partial_{\bn}^2u|_\gamma = f \quad\text{and}\quad \partial_{\bn}^3 u +(2-\nu)  \partial_{\bn} \partial^2_{\btau} u+(1-\nu) \kappa \lp\partial_{\btau}^2 u -\partial_{\bn}^2 u\rp|_\gamma = g ,
\end{equation}
where $\kappa$ is the curvature of $\gamma$, and $\btau$ is the unit tangent to $\gamma$. These boundary conditions are of particular interest in a number of applications in geophysics and structural engineering \cite{lakitosh2012analysis,kapania2000static,sergienko13,Nekrasov_MacAyeal_2023,seide1958,zhao2002plate,MEYLAND2021,Ukrainskii2018,meylan96,meylan2002,sergienko13}.

One approach for computing scattering by periodic structures relies on constructing Poincar\'e--Steklov operators, such as Robin-to-Robin \cite{turc2025} and Dirichlet-to-Neumann maps \cite{fliss2009exact,fliss2010exact,amenoagbadji2023wave,fliss2020time,fliss2019wave,klindworth2014numerical, pierre}. These operators are first computed on a \textit{unit cell} and then used to find the corresponding operator for a half-space with a semi-infinite periodic boundary. Finally, the operators are merged to solve for the field scattered by a periodic half-space. So far, these methods have been restricted to second order equations and have not been applied to models of flexural waves. This method also has the drawback that finding the half-space Poincar\'e--Steklov operator requires solving an expensive Riccati equation and can sometimes only be solved in the dissipative case.

Here, we instead use an approach based on the Floquet--Bloch transform (also known as the array scanning method \cite{munk_plane-wave_1979}). In this method, we write $u$ in terms of an inverse Floquet--Bloch transform:
\begin{equation}\label{eq:FB_tran}
    u(\bx) = \frac{d}{2\pi} \int_c u_{\xi}(\bx)\, \opd \xi,
\end{equation}
where $c$ is a complex contour connecting the points $\pm \pi/d$ and avoiding the branch cuts of $u_\xi$, defined by its \textit{quasi-periodicity} with respect to the parameter $\xi$:
\begin{equation} \label{eq:quasiperiodicity}
    u_{\xi}(\bx+d\bs e_1) = e^{i\xi d}u_{\xi}(\bx).
\end{equation}
In this work, we assume that the boundary $\gamma$ is periodic in the $\bs e_1$-direction with period $d$, the solution $u_\xi$ is quasi-periodic in the $\bs e_1$-direction, and we write $\bx = (x_1, x_2)$. Each quasi-periodic solution $u_\xi(\bx)$ can be found independently and then the total solution can be assembled using \eqref{eq:FB_tran}. This approach has been applied to a variety of Helmholtz scattering problems \cite{rana1981current,capolino2005mode,lechleiter2017convergent,agocs2024trapped,zhang2021numerical,agocs2025complex}. One popular method for solving each quasi-periodic problem is to convert it into a boundary integral equation on the restriction of $\gamma$ to the unit cell. There is a long history of such methods for Helmholtz problems \cite{agocs2024trapped,agocs2025complex,desanto1998theoretical,arens2006integral,pinto2021fast,meng2023new,bruno2014rapidly,bruno2017rapidly,perez2018domain,barnett2011new,gillman2013fast, zhang2021fast,zhang2022fast,strauszer2023windowed}.
In this work, we introduce integral equation formulations for quasi-periodic flexural wave problems based on the free-space integral equations derived in \cite{nekrasov2025boundary,lindsay2018boundary,farkas1989mathematical}

An alternate method for approaching this class of flexural problems for the special case of rectangular boundaries is given in \cite{xu2022modulated}. For rectangular boundaries it is possible to use coupled mode theory to determine solutions on each rectangular stub and match the solution to the rest of the plate. These results were validated using a finite element method with an absorbing layer away from the boundary. For periodic arrays of obstacles (not interfaces) with clamped boundary conditions, the papers \cite{li2025numerical,bao2025convergence} analyze \emph{perfectly matched layers} (PML), which can be used in conjunction with standard finite element techniques. The integral equation approach taken here works both for arbitrary smooth boundaries and obstacles, only involves discretizing one-dimensional boundaries rather than areas, leads to linear systems which are well-conditioned under refinement, and automatically enforces outgoing radiation conditions.

The remainder of this paper is structured as follows. In Section~\ref{sec:green} we introduce the quasi-periodic Green's functions. In Section~\ref{sec:ie} we introduce the integral equation formulations of the quasi-periodic problems. In Section~\ref{sec:numeric_demo} we demonstrate a method for discretizing and solving these integral equation formulations and evaluating the inverse Floquet--Bloch transform. We then test our solver on a few examples. In Section~\ref{sec:trapped} we describe how to find trapped modes and their dispersion relations. In Section~\ref{sec:merge} we extend the method of \cite{agocs2025complex} to simulating the scattering from two semi-infinite arrays of objects. Finally, in Section~\ref{sec:conclusion} we give some concluding remarks.

\section{Quasi-periodic Green's functions} \label{sec:green}
The boundary integral formulations for quasi-periodic flexural problems with clamped, supported, and free plate boundary conditions are constructed using the quasi-periodic Green's functions for the Helmholtz, flexural, and Laplace equations. In this section, we define these quasi-periodic Green's functions and provide several formulas which are useful in their evaluation.

\subsection{Helmholtz Green's function}
The (outgoing) quasi-periodic Helmholtz Green's function is the solution to the PDE
\begin{equation}
\begin{cases}
    \Delta G_{\xi,k,H}(\bx) + k^2 G_{\xi,k,H}(\bx) = \delta(\bx),\\
    G_{\xi,k,H}(\bx+ d \bs e_1) = G_{\xi,k,H}(\bx)e^{i\xi d},
\end{cases}
\end{equation}
augmented with appropriate radiation conditions as $|x_2|\to \infty$. It can be written explicitly as the following sum (equation 2.3 in~\cite{linton2010lattice})
\begin{equation}\label{eq:Ghelm_cond}
         G_{\xi,k,H}(\bx) = \sum_{n=-\infty}^{\infty} e^{in\xi d}G_{k,H}(\bx + nd\bs e_1),
\end{equation}
which is conditionally convergent for real $k.$ Here $G_{k,H}$ denotes the free-space (non-quasi-periodic) Helmholtz Green's function
\begin{equation}
    G_{k,H}(\bx):= \frac{i}{4}H_0^{(1)}(k\|\bx\|).
\end{equation}
Due to its conditional convergence for real $k$, the formula in  \eqref{eq:Ghelm_cond} is impractical for numerical purposes. Instead we compute $G_{\xi,k,H}$ in one of two ways. 
For $|x_2|$ large, it is convenient to write $G_{\xi,k,H}$ as a Bloch series, using equation 2.9 in~\cite{linton2010lattice}:
\begin{equation}
    G_{\xi,k,H} (\bx) = \sum_{m=-\infty}^\infty e^{i\xi_m x_1}\frac{e^{\alpha(\xi_m)|x_2|}}{-2\alpha(\xi_m)} \label{eq:dualsum}
\end{equation} 
with~$\xi_m = \xi + \frac{2\pi}{d}m$ and~$\alpha(\xi) = -\sqrt{i(\xi-k)}\sqrt{- i(\xi+k)}$. Throughout this work, we take the branch cut of the square root to lie along the negative real axis. This series converges exponentially in $m,$ proportional to $e^{-\frac{2\pi |m|}d|x_2|}$, which makes it convenient whenever $|x_2|$ is bounded away from zero. Practically, we use it when $|x_2|>d$. For $|x_2|\le d$, we instead compute $G_{\xi,k}$ by building a Bessel expansion for the $|n|>l$ terms in \eqref{eq:Ghelm_cond} and directly summing the remaining terms. This results in the formula
\begin{equation}
    G_{\xi,k,H}(\bx) = \sum_{j=-l}^le^{i\xi jd}G_{k,H}(\bx+jd\bs e_1) + \frac12 S_0 J_0(k\|\bx\|) +\sum_{n=1}^\infty S_n J_n(k\|\bx\|) \cos(n\opnm{arg} \bx), \label{eq:bessel_sum}
\end{equation}
for~$\|\bx\|< (l+1)d$. The Bessel series coefficients, known as  lattice coefficients, have explicit formulas as integrals (see equation 17 in \cite{yasumoto2002efficient}) and have fast algorithms for computing them \cite{denlinger2017fast}. Choosing $l=2,$ the series in \eqref{eq:bessel_sum} converges rapidly for all $\|\bx\|<\sqrt{5/4}d$, without adding too many direct interactions in the first term of \eqref{eq:bessel_sum}. For points outside of the unit cell ($|x_1|>d/2$), we compute $G_{\xi,k,H}(\bx)$ by mapping $\bx$ back to the unit cell $|x_1|<d/2$, using \eqref{eq:bessel_sum}, and explicitly using the quasi-periodicity. We refer the reader to\cite{agocs2024trapped,agocs2025complex}, and the references therein, for further details. 

\begin{remark}
Naive implementations of the derivatives of \eqref{eq:bessel_sum} may lead to catastrophic cancellations when $\|\bx\|$ is small, even though the function $G_{\xi,k}(\bx) - G_k(\bx)$ is analytic in the unit cell. In principle this issue can be alleviated by analytically canceling singular terms which arise in the derivatives, though we do not pursue this here.
\end{remark}

\subsection{Flexural Green's function}
The Green's function for the quasi-periodic flexural problem can be constructed easily from the quasi-periodic Helmholtz Green's functions (see \cite{haslinger2014symmetry}). It satisfies the PDE
\begin{equation}
\begin{cases}
    \Delta^2 G_{\xi,k,F}(\bx) - k^4 G_{\xi,k,F}(\bx) = \delta(\bx),\\
    G_{\xi,k,F}(\bx+ d \bs e_1) = G_{\xi,k,F}(\bx)e^{i\xi d},
\end{cases}
\end{equation}
together with suitable radiation conditions as $|x_2| \to \infty.$ In order to compute this Green's function, we use the identity
\begin{equation}
    G_{\xi,k,F} =\frac{1}{2k^2}\lp G_{\xi,k,H} - G_{\xi,ik,H}\rp. \label{eq:flex_gf}
\end{equation}
The two quasi-periodic Helmholtz Green's functions can be evaluated using the approach described in the previous section. With this choice of construction, when $|x_2|>d$, $G_{\xi,k,F}$ has four sets of branch cuts, extending vertically from the points $\pm k+2n\pi/d$ and $\pm i k+2n\pi/d$ for $n\in \mathbb{Z}$.  When $|x_2|<d$, the corresponding formula for $G_{\xi,k,F}$ has curved branch cuts emanating from the same points. The integration contour $c$ in \eqref{eq:FB_tran} must be chosen to avoid both sets of branch cuts. For $k \in \lp 2n\pi/d, (2n+1) \pi/d \rp$ for some integer $n$, we achieve this by choosing a $c$ that lives in the second and fourth quadrants of the complex plane. For $k$-s not in this interval, the contour $c$ should live in the first and third quadrants, and it is necessary to add the residue associated with the poles discussed in Section~\ref{sec:trapped}. In this work, we restrict our attention to the easier case of $k$. We remark in passing that the free-space Green's function for the flexural wave equation can be defined analogously to \eqref{eq:flex_gf}.

\subsection{Laplace Green's function}
In this section, we give formulas for the quasi-periodic Laplace Green's function, which solves
\begin{equation}
\begin{cases}
    \Delta G_{\xi,L}(\bx) = \delta(\bx),\\
    G_{\xi,L}(\bx+d\bs e_1) = G_{\xi,L}(\bx)e^{i \xi d}.
\end{cases}
\end{equation}
The conditions at infinity for this problem are more subtle than in the previous cases, but any Green's function satisfying the above two equations will be sufficient for our purposes. When $|x_2|$ is large enough, the most convenient way to evaluate $G_{\xi,L}$ is to evaluate \eqref{eq:dualsum} at $k=0$ (see \cite{moroz2006quasi}). For $|x_2|$ small, we use standard techniques to build a spherical harmonic expansion (analogous to \eqref{eq:bessel_sum}). To derive this, we start by writing $G_{\xi,L}$ as the conditionally convergent sum of periodic copies of the Laplace Green's function:
\begin{equation}
     G_{\xi,L}(\bx) = \sum_{n=-\infty}^{\infty} e^{in\xi d}\frac{1}{2\pi}\log\|\bx + nd\bs e_1\|.
\end{equation}
Each term can be written as a spherical harmonic expansion \cite{joslin1983multipole}:
\begin{equation}
    \frac{1}{2\pi}\log\|\bx + nd\bs e_1\| = \frac{1}{2\pi}\log (nd) - \frac{1}{2\pi}\sum_{l=1}^\infty \frac{r^l}{l (nd)^l}\cos (l\theta),
\end{equation}
where $\theta := \arg(\bx).$ Summing over $n>0$ gives
\begin{equation}
    \sum_{n=1}^{\infty} e^{in\xi d}\frac{1}{2\pi}\log\|\bx + nd\bs e_1\| = \frac{1}{2\pi}\sum_{n=1}^{\infty} e^{in\xi d}\log (nd) - \frac{1}{2\pi}\sum_{n=1}^{\infty} e^{in\xi d}\sum_{l=1}^\infty \frac{r^l}{l (nd)^l}\cos (l\theta).
\end{equation}
If $\Im \xi<0$, then both sums converge exponentially, and we may reorder the sums as
\begin{equation}
\begin{aligned}
    \sum_{n=1}^{\infty} e^{in\xi d}\frac{1}{2\pi}\log\|\bx + nd\bs e_1\| &= \frac{1}{2\pi}\sum_{n=1}^{\infty} e^{in\xi d}\log (nd) - \frac{1}{2\pi}\sum_{l=1}^\infty \frac{r^l}{l {d}^l}\cos (l\theta)\sum_{n=1}^{\infty} \frac{e^{in\xi d}}{n^l}\\
    &=C_{\xi,+} -\frac{1}{2\pi}\sum_{l=1}^\infty \frac{r^l}{l {d}^l}\cos (l\theta) \opnm{Li}_l(e^{i\xi d}),
\end{aligned}
\end{equation}
where $C_{\xi,+}$ is independent of $\bx$ and $\opnm{Li}_l$ is the polylogarithm function. Since the right hand side is analytic, we can extend this formula to $\Im \xi=0$. Repeating the same argument for the $n<0$ terms gives the formula
\begin{equation}\label{eq:Lap_multi}
    G_{\xi,L}(\bx) = \frac{1}{2\pi}\log\|\bx\| + C_\xi - \frac{1}{2\pi}\sum_{l=1}^\infty \frac{r^l}{l {d}^l}\cos (l\theta) \lp \opnm{Li}_l(e^{i\xi d})+(-1)^l\opnm{Li}_l(e^{-i\xi d})\rp.
\end{equation}

Rather than derive an explicit analytic expression for the constant $C_\xi,$ in our calculations we calculate it numerically by comparing \eqref{eq:Lap_multi} with the Laplace analog of \eqref{eq:dualsum} at a point in which both expansions are rapidly convergent. In practice, we add and subtract a few periodic copies of the point source to ensure that the spherical harmonic expansion converges rapidly (analogous to \eqref{eq:bessel_sum}).

The boundary integral formulations of the free-space flexural problem developed in \cite{nekrasov2025boundary} use the Poincar\'e--Bertrand formula \cite{muskhelishvili2008singular,hang2009generalized} to simplify the resulting integral equations. It will therefore be useful to have a quasi-periodic analog when deriving our integral equations below. We propose the following conjecture, which gives a quasi-periodic version of the standard identity.

\begin{conjecture}[Quasi-periodic Poincar\'e--Bertrand formula]
    Let $\gammaunit$ denote the portion of $\gamma$ lying in the strip $[-d/2,d/2) \times \mathbb{R}$ and  
    \begin{equation}
        \mathcal{H}_\xi[\sigma](\bx) = \int_\gammaunit \partial_{\btau(\by)} G_{\xi,L}(\bx-\by)\sigma(\by)\,{\rm{d}}\by \quad \text{and}\quad \mathcal{D}_\xi[\sigma](\bx) = \int_\gammaunit \partial_{\bn(\by)} G_{\xi,L}(\bx-\by)\sigma(\by)\,{\rm{d}}\by ,
    \end{equation}
    be the quasi-periodic Hilbert transform and Laplace double layer potential respectively. If $\sigma$ is H\"older continuous and quasi-periodic, then
    \begin{equation}
        \frac14\mathcal{H}_\xi^2[\sigma] = -\frac{\sigma}4 + \mathcal{D}_\xi^2[\sigma].
    \end{equation}
\end{conjecture}
This formula in its non-quasiperiodic form is well-known to hold for a closed simple smooth curve \cite{hang2009generalized}. In numerical experiments, over a number of examples we observe that the conjectured identity holds with a relative accuracy of at least $1.7\times 10^{-8}$ when applied to a smooth quasi-periodic right hand side.

\section{Integral equations for quasi-periodic flexural problems}\label{sec:ie}
In this section, we describe our method for solving quasi-periodic flexural boundary value problems with clamped, supported, and free plate boundary conditions. As discussed in the introduction, we solve the quasi-periodic problem separately for each Bloch parameter $\xi$ and then write the total solution $u$ using the inverse Floquet--Bloch transform \eqref{eq:FB_tran}. Recall that the solution $u_\xi$ to the quasi-periodic flexural problem satisfies the following equation: 
\begin{equation}\label{eq:quaspde}
\begin{cases}
    \Delta^2 u_\xi - k^4 u_\xi= 0 , \\
    u_{\xi}(\bx+d\bs e_1) = e^{i\xi d}u_{\xi}(\bx),
\end{cases}
\end{equation}
paired with the appropriate boundary conditions. In the remainder of this section, we give integral equation formulations for a few different boundary conditions. All of our representations are obtained by inserting quasi-periodic kernels into the representations introduced in \cite{nekrasov2025boundary,lindsay2018boundary,farkas1989mathematical}. We include these integral representations here for completeness. To avoid notational clutter we suppress the dependence of $G_{\xi,k,F}$ on $k,$ writing $G_\xi$ instead. Note that the integral formulations presented here require that the boundary $\gamma$ possesses sufficient regularity. Throughout this section, $\bx$ and $\by$ denote points inside the unit cell. The solution outside of the unit cell may be obtained via the quasi-periodicity condition.
We also let $f_\xi$ and $g_\xi$ be the Floquet--Bloch transform of the boundary data $f$ and $g$ defined in Section~\ref{sec:intro}.

\subsection{Clamped plate boundary conditions}
As stated in the introduction, the clamped plate boundary conditions are
\begin{equation}
    u_\xi|_\gammaunit = f_\xi \quad\text{and}\quad \partial_{\bn} u_\xi|_\gammaunit = g_\xi . \, \label{eq:clamped_quasi}
\end{equation}
Following \cite{nekrasov2025boundary,lindsay2018boundary,farkas1989mathematical}, we suppose that $u_\xi$ can be written as
\begin{multline}\label{eq:clamp_rep}
   u_\xi(\bx) = \int_\gammaunit \left[\lp \partial_{\bn(\by)}^3 G_\xi(\bx-\by)+ 3\partial_{\bn(\by)}\partial_{\btau(\by)}^2 G_\xi(\bx-\by) \rp  \sigma(\by) \right. \\+ \left. \lp -\partial_{\bn(\by)}^2 G_\xi(\bx-\by)+ \partial_{\btau(\by)}^2 G_\xi(\bx-\by) \rp  \rho(\by)\right]\,{\rm{d}} \by ,
\end{multline}
for some unknown densities $\sigma,\rho$. By construction, $u_\xi$ will satisfy the quasi-periodic flexural equation for any $\sigma,\rho$. To enforce the boundary conditions \eqref{eq:clamped_quasi}, we note that \eqref{eq:bessel_sum} implies that $G_\xi$ has the same singularity as $\|\bx-\by\|\to 0$ as the flexural Green's function. The quasi-periodic layer potentials therefore satisfy the same jump relations as those based on the free-space flexural Green's function. Thus, we find that $u_\xi$ given by \eqref{eq:clamp_rep} satisfies the clamped plate boundary conditions provided that $\sigma$ and $\rho$ satisfy
\begin{equation}\label{eq:clamp_IE}
    \begin{pmatrix}
        -\frac12 I & 0 \\
        \kappa I & - \frac12 I
    \end{pmatrix}\begin{pmatrix}
        \sigma\\\rho
    \end{pmatrix} + \begin{pmatrix}
        \mathcal{K}_{11} & \mathcal{K}_{12} \\
        \mathcal{K}_{12} & \mathcal{K}_{22}
    \end{pmatrix} \begin{pmatrix}
        \sigma\\\rho
    \end{pmatrix} = \begin{pmatrix}
        f_\xi\\g_\xi
    \end{pmatrix},
\end{equation}
where $\kappa$ is the curvature of $\gamma$ and $\mathcal{K}_{11},\mathcal{K}_{12},\mathcal{K}_{21},\mathcal{K}_{22}$ are integral operators with kernels
\begin{equation}
\begin{aligned}
  K_{11}(\bx,\by) &= \partial_{\bn(\by)}^3G_\xi(\bx-\by)+ 3\partial_{\bn(\by)}\partial_{\btau(\by)}^2 G_\xi(\bx-\by) ,  \\
  K_{12}(\bx,\by) &=-\partial_{\bn(\by)}^2 G_\xi(\bx-\by)+ \partial_{\btau(\by)}^2 G_\xi(\bx-\by) ,\\
  K_{21}(\bx,\by) &=\partial_{\bn(\bx)}\partial_{\bn(\by)}^3G_\xi(\bx-\by)+ 3\partial_{\bn(\bx)}\partial_{\bn(\by)}\partial_{\btau(\by)}^2 G_\xi(\bx-\by) ,\\
  K_{22}(\bx,\by) &=-\partial_{\bn(\bx)}\partial_{\bn(\by)}^2 G_\xi(\bx-\by)+ \partial_{\bn(\bx)}\partial_{\btau(\by)}^2 G_\xi(\bx-\by) ,
\end{aligned}
\end{equation}
respectively. It was observed in \cite{nekrasov2025boundary} that the kernels corresponding to the free-space Green's function are all continuous on surface for a smooth boundary. Since the quasi-periodic kernels are a smooth correction to the free-space kernels, they are all sufficiently smooth in our case. Equation~\eqref{eq:clamp_IE} is thus a Fredholm second kind integral equation.

\subsection{Supported plate boundary conditions}

As stated in the introduction, the supported plate boundary conditions are given by
\begin{equation*}
    u_\xi|_\gammaunit = f_\xi \quad\text{and}\quad \nu \Delta u_\xi +(1-\nu) \partial_{\bn}^2u_\xi |_\gammaunit = g_\xi \, ,
\end{equation*}
Following \cite{nekrasov2025boundary}, we write $u_\xi$ as
\begin{multline}\label{eq:supported_rep}
    u_\xi (\bx) = \int_\gammaunit \bigg[ \big( \partial_{\bn(\by)}^3 G_\xi (\bx-\by) + \alpha_1 \partial_{\bn(\by)} \partial_{\btau(\by)}^2 G_\xi (\bx-\by) + \alpha_2 \kappa(\by) \partial_{\bn(\by)}^2 G_\xi(\bx-\by) \\
    + \alpha_3 \kappa'(\by) \partial_{\btau(\by)} G_\xi(\bx-\by) \big) \sigma(\by) 
    + \partial_{\bn(\by)} G_\xi(\bx-\by) \rho(\by) \bigg] \, \dd \by \, ,
\end{multline}
where $\kappa'(\by)$ denotes the arclength derivative of curvature at the point $\by \in \gammaunit$ and the coefficients $\alpha_1,\alpha_2,$ and $\alpha_3$ are given by 
\begin{align*}
    \alpha_1 = 2- \nu \, ,  \ \ \alpha_2 = \frac{(-1+\nu)(7+\nu)}{3-\nu} \, , \ \ \alpha_3 = \frac{(1-\nu)(3+\nu)}{1+\nu} \, .
\end{align*} These layer potentials also satisfy the same jump relations as their free-space analogs. Thus $u_\xi$ satisfies the supported plate boundary conditions provided $\sigma$ and $\tau$ satisfy the boundary integral equation
\begin{equation}\label{eq:supported_IE}
    \begin{pmatrix}
        -\frac12 I & 0 \\
        c_0 \kappa^2 I & - \frac12 I
    \end{pmatrix}\begin{pmatrix}
        \sigma\\\rho
    \end{pmatrix} + \begin{pmatrix}
        \mathcal{K}_{11} & \mathcal{K}_{12} \\
        \mathcal{K}_{21} & \mathcal{K}_{22}
    \end{pmatrix} \begin{pmatrix}
        \sigma\\\rho
    \end{pmatrix} = \begin{pmatrix}
        f_\xi\\g_\xi
    \end{pmatrix},
\end{equation}
where $c_0 $ is given by
\begin{equation*}
    c_0 =  \frac{(\nu-1)(\nu+3)(2\nu-1)}{2(3-\nu)} \, ,
\end{equation*}
and $\cK_{11},\cK_{12},\cK_{21},$ and $\cK_{22}$ are boundary integral operators with kernels
\begin{align*}
    K_{11}(\bx,\by) &=  \left[\partial_{\bn(\by)}^3 + \alpha_1 \partial_{\bn(\by)} \partial_{\btau(\by)}^2  + \alpha_2 \kappa(\by) \partial_{\bn(\by)}^2
    + \alpha_3 \kappa'(\by) \partial_{\btau(\by)}\right] G_\xi(\bx-\by) \, , \\
    K_{12}(\bx,\by) &= \partial_{\bn(\by)} G_\xi(\bx-\by) \, , \\
    K_{21}(\bx,\by) &=  \left[\partial_{\bn(\bx)}^2 \partial_{\bn(\by)}^3  + \alpha_1 \partial_{\bn(\bx)}^2 \partial_{\bn(\by)} \partial_{\btau(\by)}^2 + \alpha_2 \kappa(\by) \partial_{\bn(\bx)}^2 \partial_{\bn(\by)}^2  \right.   + \alpha_3 \kappa'(\by) \partial_{\bn(\bx)}^2 \partial_{\btau(\by)}   \\
    &\quad  +\nu \partial_{\btau(\bx)}^2 \partial_{\bn(\by)}^3 +  \nu \alpha_1 \partial_{\btau(\bx)}^2 \partial_{\bn(\by)} \partial_{\btau(\by)}^2  + \nu \alpha_2 \kappa(\by) \partial_{\btau(\bx)}^2 \partial_{\bn(\by)}^2
    \left.+ \nu \alpha_3 \kappa'(\by) \partial_{\btau(\bx)}^2 \partial_{\btau(\by)} \right]G_\xi(\bx-\by) \, , \\
    K_{22}(\bx,\by) &= \left[\partial_{\bn(\bx)}^2 \partial_{\bn(\by)}+ \nu \partial_{\btau(\bx)}^2 \partial_{\bn(\by)}\right] G_\xi(\bx-\by) \, , 
\end{align*}
respectively. Since the free-space kernels were shown to be continuous in \cite{nekrasov2025boundary}, and the quasi-periodic kernels are given by the free-space kernels plus a smooth correction, the integral operators in \eqref{eq:supported_IE} are compact and the equation is Fredholm second kind. 

\subsection{Free plate boundary conditions}
As stated in the introduction, the free plate boundary conditions are given by
\begin{equation*}
    \nu \Delta u_\xi +(1-\nu) \partial_{\bn}^2u_\xi|_\gammaunit = f_\xi \quad\text{and}\quad \partial_{\bn}^3 u_\xi +(2-\nu)  \partial_{\bn} \partial^2_{\btau} u_\xi+(1-\nu) \kappa \lp\partial_{\btau}^2 u_\xi -\partial_{\bn}^2 u_\xi\rp|_\gammaunit = g_\xi ,
\end{equation*}
on $\gammaunit$.
To solve the free plate problem, we follow \cite{nekrasov2025boundary} and suppose that $u_\xi$ can be written as
\begin{equation}\label{eq:free_rep}
   u_\xi(\bx) = \int_\gammaunit \left[ \partial_{\bn(\by)}G_\xi(\bx-\by)  \sigma(\by) +\frac{1+\nu}2 \partial_{\btau(\by)}G_\xi(\bx-\by)  \mathcal{H}_\xi[\sigma](\by) +G_\xi(\bx-\by)   \rho(\by)\right]\,{\rm{d}} \by,
\end{equation}
which satisfies the quasi-periodic flexural wave equation in $\Omega$ for any~$\sigma,\tau$. It is easily seen that these quasi-periodic layer potentials also satisfy the same jump relations as their free-space analogs. Thus $u_\xi$ satisfies the free plate boundary conditions provided $\sigma$ and $\tau$ satisfy
\begin{equation}\label{eq:Free_IE1}
    \begin{pmatrix}
        -\frac{1}2I & 0\\
        -\frac{1+\nu}4 \partial_\tau \mathcal{H}_\xi & \frac{1}2I
    \end{pmatrix}\begin{pmatrix}
        \sigma\\ \tau
    \end{pmatrix} + \begin{pmatrix}
        \mathcal{K}_{11}^a + \frac{1+\nu}2\mathcal{K}_{11}^b\mathcal{H}_\xi & \mathcal{K}_{12} \\
        \mathcal{K}_{21}^a + \frac{1+\nu}2\mathcal{K}_{21}^b\mathcal{H}_\xi & \mathcal{K}_{22}
    \end{pmatrix} \begin{pmatrix}
        \sigma\\\rho
    \end{pmatrix} = \begin{pmatrix}
        f_\xi\\g_\xi
    \end{pmatrix},
\end{equation}
where $\mathcal{K}_{11}^a,\mathcal{K}_{11}^b,\mathcal{K}_{12},\mathcal{K}_{21}^a,\mathcal{K}_{21}^b,\mathcal{K}_{22}$ are integral operators with kernels
\begin{equation*}
\begin{aligned}
  K_{11}^a(\bx,\by) &= \left[\partial_{\bn(\bx)}^2\partial_{\bn(\by)}+\nu\partial_{\btau(\bx)}^2 \partial_{\bn(\by)}\right]G_\xi(\bx-\by) ,  \\
    K_{11}^b(\bx,\by) &= \left[\partial_{\bn(\bx)}^2\partial_{\btau(\by)}+\nu\partial_{\btau(\bx)}^2 \partial_{\btau(\by)}\right]G_\xi(\bx-\by) ,  \\
  K_{12}(\bx,\by) &=\left[\partial_{\bn(\bx)}^2+\nu\partial_{\btau(\bx)}^2 \right]G_\xi(\bx-\by)  ,\\
  K_{21}^a(\bx,\by) &=\left[\partial_{\bn(\bx)}^3\partial_{\bn(\by)}+(2-\nu)\partial_{\bn(\bx)}\partial_{\btau(\bx)}^2\partial_{\bn(\by)}+(1-\nu)\kappa(\bx)\lp \partial_{\btau(\bx)}^2 \partial_{\bn(\by)}- \partial_{\bn(\bx)}^2\partial_{\bn(\by)}\rp\right] G_\xi(\bx-\by),\\
  K_{21}^b(\bx,\by) &=\left[\partial_{\bn(\bx)}^3\partial_{\btau(\by)}+(2-\nu)\partial_{\bn(\bx)}\partial_{\btau(\bx)}^2\partial_{\btau(\by)}+(1-\nu)\kappa(\bx)\lp \partial_{\btau(\bx)}^2 \partial_{\btau(\by)}- \partial_{\bn(\bx)}^2\partial_{\btau(\by)}\rp\right] G_\xi(\bx-\by),\\
  K_{22}(\bx,\by) &=\left[\partial_{\bn(\bx)}^3+(2-\nu)\partial_{\bn(\bx)}\partial_{\btau(\bx)}^2+(1-\nu)\kappa(\bx)\lp \partial_{\btau(\bx)}^2- \partial_{\bn(\bx)}^2\rp\right] G_\xi(\bx-\by),
\end{aligned}
\end{equation*}
respectively. If we apply the quasi-periodic Poincar\'e--Bertrand formula, then \eqref{eq:Free_IE1} can be written as
\begin{equation}\label{eq:Free_IE2}
    \begin{pmatrix}
        \lp-\frac{1}2+\frac{(1+\nu)^2}{8}\rp I & 0\\
        0 & \frac{1}2I
    \end{pmatrix}\begin{pmatrix}
        \sigma\\ \tau
    \end{pmatrix} + \begin{pmatrix}
        \mathcal{K}_{11}^a + \frac{1+\nu}2\lp \mathcal{K}_{11}^b+\frac{1+\nu}4\mathcal{H}_\xi\rp\mathcal{H}_\xi - \frac{(1+\nu)^2}4\mathcal{D}_\xi & \mathcal{K}_{12} \\
        -\frac{1+\nu}4 \partial_\tau \mathcal{H}_\xi+\mathcal{K}_{21}^a + \frac{1+\nu}2\mathcal{K}_{21}^b\mathcal{H}_\xi & \mathcal{K}_{22}
    \end{pmatrix} \begin{pmatrix}
        \sigma\\\rho
    \end{pmatrix} = \begin{pmatrix}
        f_\xi\\g_\xi
    \end{pmatrix}.
\end{equation}
Since the quasi-periodic kernels are given by the free-space kernels plus a correction, each of the above integral operators are free-space operator plus a correction that is compact on $L^2(\gammaunit)$. Since all of the operators are bounded on the same space, we have that \eqref{eq:Free_IE2} is Fredholm second kind on $L^2(\gammaunit) \times L^2(\gammaunit)$.

\section{Numerical experiments}\label{sec:numeric_demo}
To evaluate the inverse Floquet--Bloch transform, we choose the contour $c$ to be given by $\Im \xi = -0.3\sin \lp d\Re \xi \rp$. The contour $c$ in \eqref{eq:FB_tran} is then discretized using a 61 point periodic trapezoid rule. Using an odd number of nodes guarantees that we avoid the removable singularity in the free plate equations at $\xi=0$. We discretize the integral equations~\eqref{eq:clamp_IE}, \eqref{eq:supported_IE}, and \eqref{eq:Free_IE2} at each $\xi$ using chunkIE \cite{askham_chunkie_2024}. This package splits $\gammaunit$ into 16th order Gauss-Legendre panels and discretizes the integral operators using a collocation scheme. To evaluate the corresponding representations \eqref{eq:clamp_rep}, \eqref{eq:supported_rep}, and \eqref{eq:free_rep}, we split the kennels into the singular free-space part and the smooth correction. We integrate the smooth part with standard Gauss-Legendre quadrature. For the singular part, we use adaptive integration. Since the singular part is the same for every $\xi$, we can reuse the adaptive integration giving us an efficient scheme. Further details on the discretization of this flavor of method are given in \cite{agocs2024trapped,agocs2025complex}.

\begin{figure}
    \centering
    \begin{subfigure}[b]{0.4\textwidth}
    \centering
        \includegraphics[width=\textwidth]{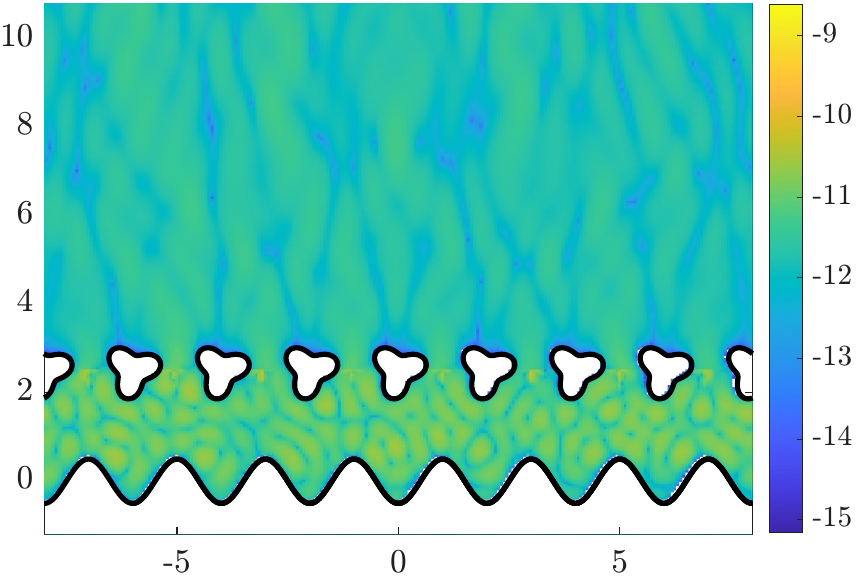}
        \caption{Clamped plate boundary conditions}
    \end{subfigure}%
    \hspace{0.5 cm}
    \begin{subfigure}[b]{0.4\textwidth}
    \centering
        \includegraphics[width=\textwidth]{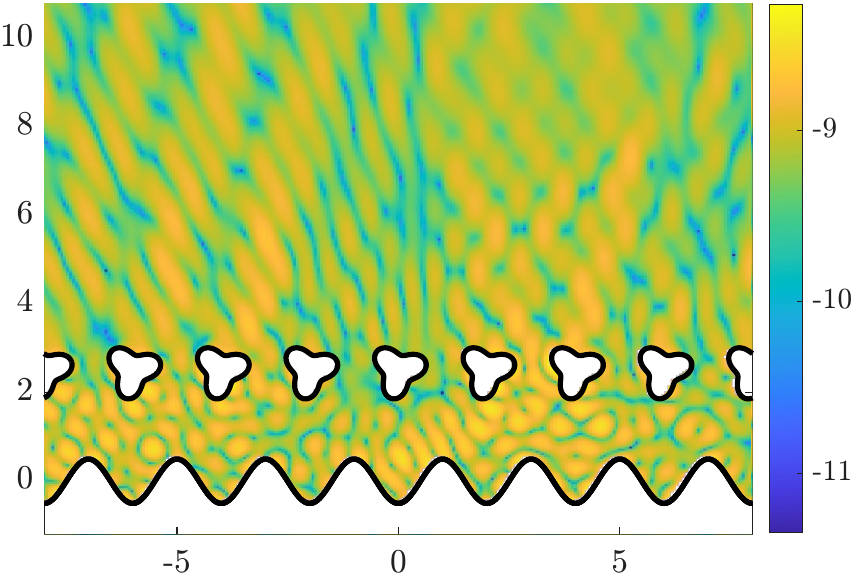}
        \caption{Free plate boundary conditions}
    \end{subfigure}
    \par\medskip 
    \begin{subfigure}[b]{0.4\textwidth}
    \centering
        \includegraphics[width=\textwidth]{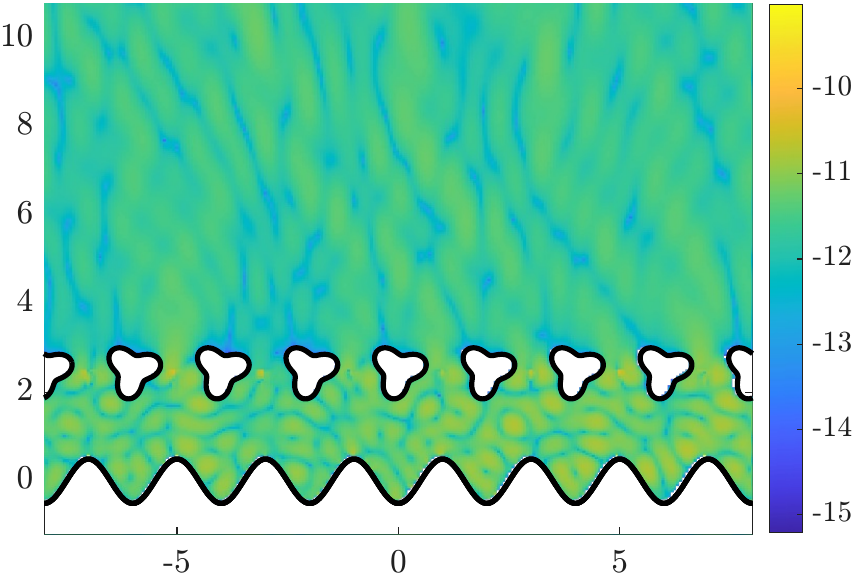}
        \caption{Supported plate boundary conditions}
    \end{subfigure}%
    \caption{These figures show $\log_{10} |u+u^{\mathrm{in}}|/\max|u^{\mathrm{in}}|$ for $\bx_0=(0,-2)$ with $k=7,$ and $d=2$.}
    \label{fig:acc_test} 
\end{figure}

To test the accuracy of our method, we let $u^{\mathrm{in}}(\bx) = G(\bx-\bx_0)$ be the free-space solution due to a point source located at $\bx=(0,-2),$ \emph{below} the interface. We then choose our boundary data for $u$ so that when $u+u^{\mathrm{in}}$ is inserted into the boundary conditions, the result is zero. Assuming uniqueness of solutions, it is then easily seen that $u=-u^{\mathrm{in}}.$ The resulting errors are shown in Figure~\ref{fig:acc_test}.

\begin{figure}
    \centering
    \begin{subfigure}[b]{0.4\textwidth}
    \centering
        \includegraphics[width=\textwidth]{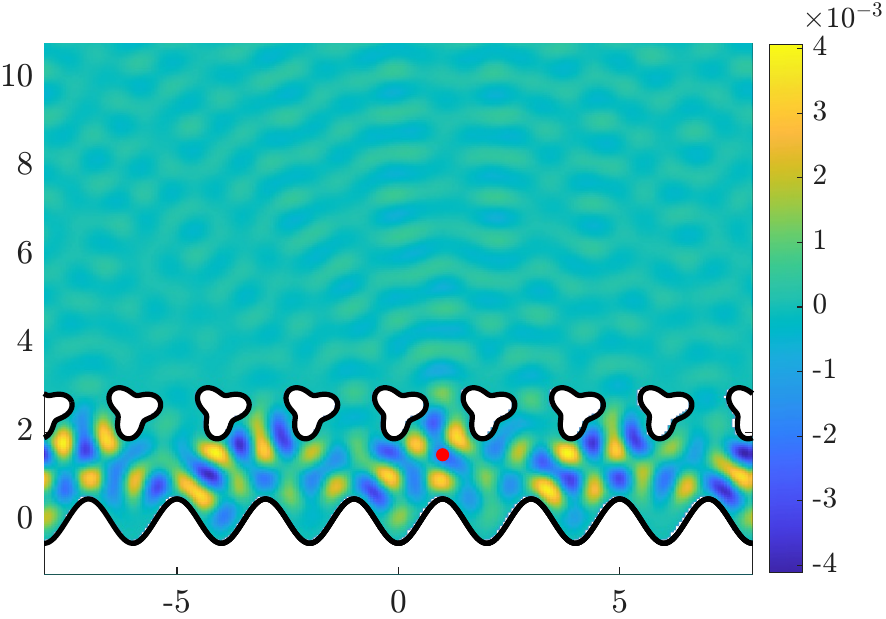}
        \caption{Clamped plate boundary conditions}
    \end{subfigure}%
    \hspace{0.5 cm}
    \begin{subfigure}[b]{0.4\textwidth}
    \centering
        \includegraphics[width=\textwidth]{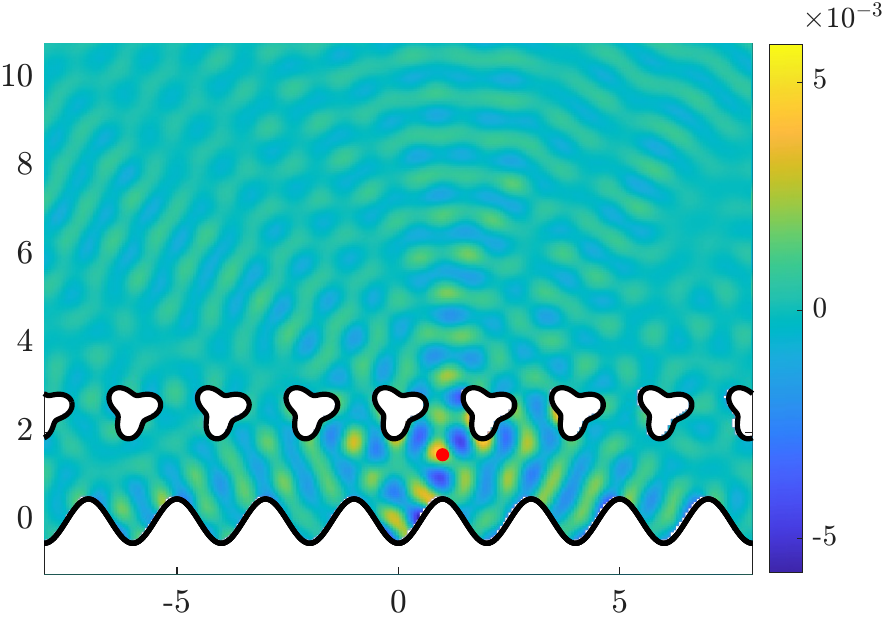}
        \caption{Free plate boundary conditions}
    \end{subfigure}
    \par\medskip
    \begin{subfigure}[b]{0.4\textwidth}
    \centering
        \includegraphics[width=\textwidth]{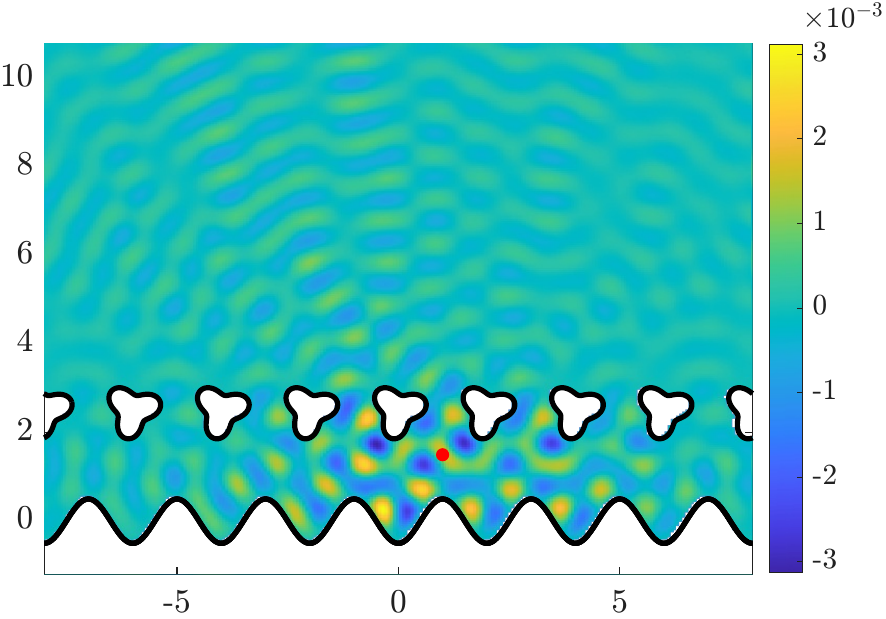}
        \caption{Supported plate boundary conditions}
    \end{subfigure}
    \caption{The imaginary parts of fields generated by a point source at $(1,1.5)$ with $k=7,$ and $d=2$.}
    \label{fig:examp_sol}
\end{figure}

 As an illustration of the performance of the method on a more physical scattering problem, we next move the source to $\bx_0=(1,1.5),$ above the grating, and again choose the boundary data so that $u+u^{\mathrm{in}}$ solves the problem with homogeneous boundary conditions. The fields due to a point source at $\bx_0=(1,1.5)$ are shown in Figure~\ref{fig:examp_sol}.

\section{Edge modes}\label{sec:trapped}

Of particular interest in a number of applications is the existence of waves that propagate along the grating over infinitely long distances without decay \cite{konenkov1960rayleigh,sinha1974some,evans2008flexural,wang2018topological,mousavi2015topologically,miniaci2018experimental}, known as `edge modes' or `trapped modes'. Upon applying the Floquet--Bloch transform they manifest as a nullspace of the quasi-periodic boundary value problems at discrete values of the parameter $\xi,$ depending on the type of boundary condition, the geometry, and the wavenumber. These nullspaces of the PDE at certain $\xi$ in turn correspond to nullspaces of the integral equation (\ref{eq:clamp_IE}), (\ref{eq:supported_IE}), or (\ref{eq:Free_IE2}) at the same $\xi$. We note that while each trapped mode corresponds to a one-dimensional nullspace at a particular $\xi$, the converse is not necessarily true. That is, it is possible that the integral equations can exhibit spurious resonances.

To find pairs $(\xi,k)$ such that the integral equations \eqref{eq:clamp_IE}, \eqref{eq:supported_IE}, or \eqref{eq:Free_IE2} exhibit nullspaces, we follow \cite{zhao2015robust} and look for zeros of the determinant of the discretized system. To find these zeros, we build an adaptive piecewise Chebyshev series for the determinant as a function of $k$ at fixed $\xi$. We then then find the roots of each  expansion by finding the roots of the associated \emph{colleague matrix} \cite{good1961colleague}. 

\begin{figure}
    \centering
    \includegraphics[width=0.6\linewidth]{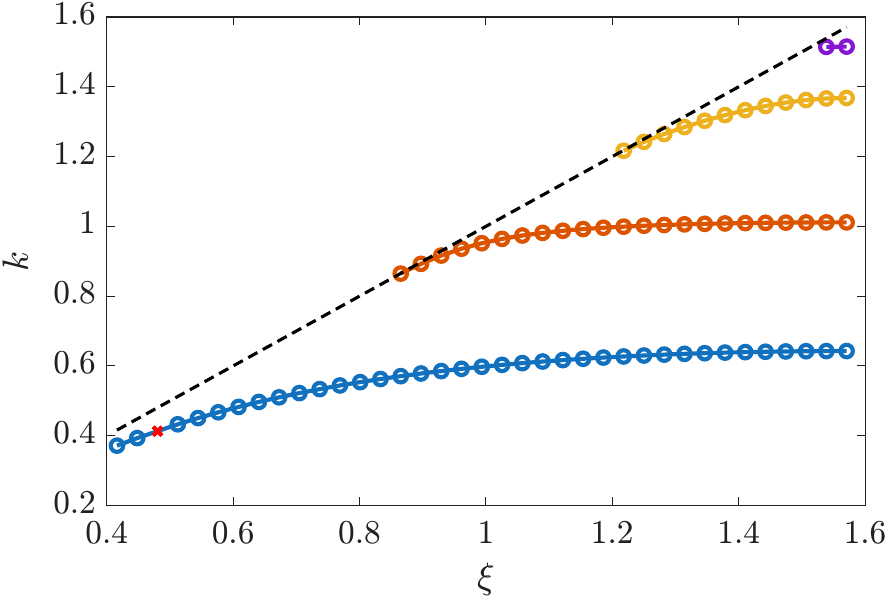}
    \caption{The dispersion relations for the free plate modes in the region shown in \figref{fig:free_modes}. The red x indicates a $\xi$ where the solver failed to find a mode with adaptive refinement tolerance $10^{-4}$. The location was obtained by linear interpolation from the adjacent $k$-s. The ratio of the first and last singular value at that $(\xi,k)$ pair was $1.9\times 10^{-7}$, which is only slightly larger than the comparable to the ratio of $1.3\times 10^{-8}$ at the adjacent pair on the left and indicates the presence of a nearby mode.}
    \label{fig:free_disp}
\end{figure}

\begin{figure}
    \centering
    \begin{subfigure}[t]{0.49\textwidth}
    \centering
    \includegraphics[width=\linewidth]{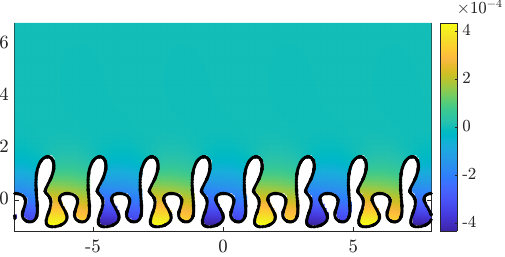}
    \end{subfigure}\hfill
       \begin{subfigure}[t]{0.49\textwidth}
    \centering
    \includegraphics[width=\linewidth]{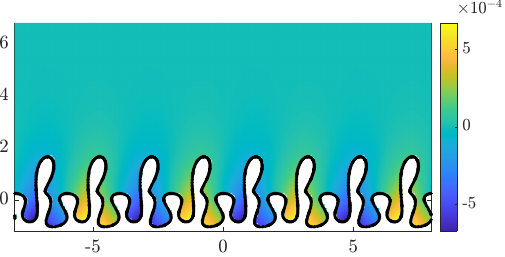}
    \end{subfigure}
    \par \medskip
       \begin{subfigure}[t]{0.49\textwidth}
    \centering
    \includegraphics[width=\linewidth]{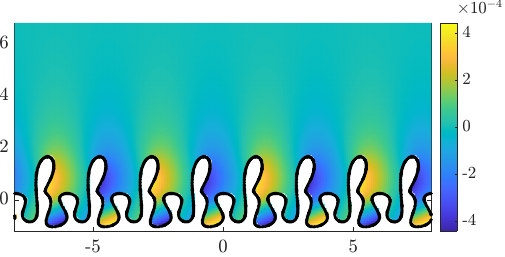}
    \end{subfigure}\hfill
       \begin{subfigure}[t]{0.49\textwidth}
    \centering
    \includegraphics[width=\linewidth]{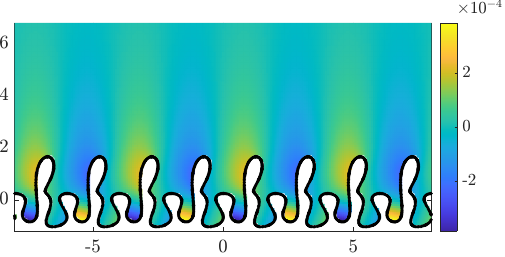}
    \end{subfigure}

    \caption{This figure shows the imaginary part of the trapped modes with $\xi=\pi/d$ and $k\in(0.3,\pi/d)$ when free plate boundary conditions are imposed on the shown $\gamma$. The frequencies are $k =$  0.642700 (top left), 1.011460  (top right), 1.367986 (bottom left), and  1.515168 (bottom right). }
    \label{fig:free_modes}
\end{figure}

\begin{figure}
\centering
\begin{subfigure}[t]{0.49\textwidth}
    \centering
    \includegraphics[width=\linewidth]{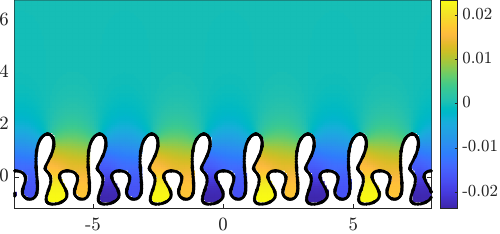}
\end{subfigure}
\hfill
\begin{subfigure}[t]{0.49\textwidth}
    \centering
    \includegraphics[width=\linewidth]{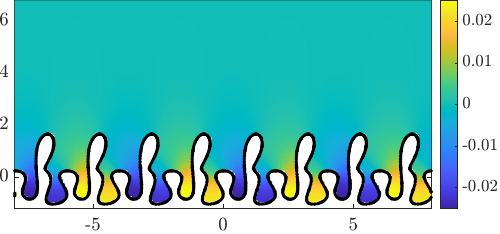}
\end{subfigure}
\par \medskip
\begin{subfigure}[t]{0.49\textwidth}
    \centering
    \includegraphics[width=\linewidth]{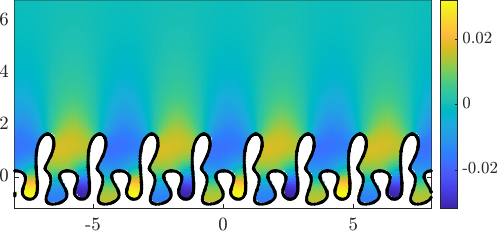}
\end{subfigure}
    \caption{The figure shows the imaginary part of the trapped modes for the Helmholtz equation with Neumann boundary conditions with $\xi = \pi/d$ and $k\in(0.3,\pi/d)$. The frequencies used are $k = 0.533608$ (left),  $k=0.9136840$ (right), and $k=1.383154$ (bottom).}
    \label{fig:neu_modes}
\end{figure}

 For the free plate problem, we employ the above method to compute the relationship between $\xi$ and $k$ for the boundary in Figure~\ref{fig:free_modes}. This dispersion relation is shown in Figure~\ref{fig:free_disp}. As our formulas for the quasi-periodic Green's functions blow up when $\xi=k$, we only look for modes with $|\xi-k|>10^{-4}$. We also do not look for modes with $k<0.38$. The modes at $\xi=\pi/d$ are shown in Figure~\ref{fig:neu_modes}. For comparison, in Figure~\ref{fig:neu_modes} we also show the modes for the Helmholtz equation with Neumann boundary conditions. We can see that the lowest frequency modes are similar to the free plate modes, but the higher frequency ones are quite different.

When the same method is applied to \eqref{eq:clamp_IE} and \eqref{eq:supported_IE}, we also find null spaces for various $k<\pi/d$ for the boundaries in \figref{fig:examp_sol} and \figref{fig:free_modes}. These corresponding null vectors generate zero field, so do not represent trapped modes. They also do not affect the accuracy of the total aperiodic solver, provided they are not too close to the contour $c$.

%

\section{Semi-infinite arrays of obstacles}\label{sec:merge} 

There has recently been great interest in periodic structures or tilings with changes in the obstacles and their effect on wave scattering (see \cite{agocs2025complex,turc2025} and the references therein). In this section, we use the integral formulations above and the complex scaling method introduced in \cite{agocs2025complex} to simulate scattering for semi-infinite arrays of compact obstacles of different types. More technical details are given in \cite{agocs2025complex} and many of the contained proofs extend \textit{mutatis mutandis} to this setting.

\subsection{Matched integral equation}
We suppose that $\gamma_{L,R}$ represent two semi-infinite periodic arrays of compact obstacles to the left and right of the line $\Gamma=\{x_1=0\}$, respectively. We let $\Omega_{L,R}\subset \bbR^2$ be the regions to the left and right of $\Gamma$ exterior to $\gamma_{L,R}$. We seek solutions $u_{L,R}$ to the following PDEs
\begin{equation}\label{eq:tot_PDE}
    \begin{cases}
        \Delta^2 u_{L} -k^4 u_{L} = 0 , \ & \bx \in \Omega_{L}\,,\\
        \Delta^2 u_{R} - k^4 u_{R} = 0 , \ & \bx \in \Omega_{R}\,,\\
    \end{cases}
\end{equation}
satisfying the continuity conditions on the fictitious interface $\Gamma$,
\begin{equation}\label{eq:cont}
    \begin{cases}
    u_R - u_L = r_0 , &\text{on } \Gamma\, ,\\
        \partial_{x_1}u_R - \partial_{x_1}u_L = r_1 , &\text{on } \Gamma\,, \\
        \partial_{x_1}^2u_R - \partial_{x_1}^2u_L = r_2  ,&\text{on } \Gamma\,,\\
        \partial_{x_1}^3u_R - \partial_{x_1}^3u_L = r_3 ,&\text{on } \Gamma\,,
    \end{cases}
\end{equation}
together with clamped, supported, or free plate boundary conditions (or a combination) on $\gamma_{L,R}$. The data $r_0,\ldots, r_3$ represent incoming data, e.g. the field due to a point source in the left or right half and its derivatives.

We proceed by reducing \eqref{eq:tot_PDE} and \eqref{eq:cont} to an integral equation on $\Gamma$. First, we let $G_{L,R}$ denote the fundamental solutions for the left and right regions satisfying the boundary conditions on $\gamma_{L,R}$. It is convenient to split
\begin{equation}
    G_{L,R}(\bx,\by) = G(\bx-\by) + w_{L,R}(\bx,\by),
\end{equation}
where $G$ is the free-space flexural Green's function, defined equivalently to \eqref{eq:flex_gf}, and $w_{L,R}$ is the scattered field, which can be computed using the method described above with boundary data $f_\xi,g_\xi$ given by plugging $G_\xi$ into the boundary conditions. We then represent $u_{L,R}$ using layer potentials with kernels constructed from these fundamental solutions
\begin{equation}\label{eq:matchrep}
    u_{L,R}(\bx) =\sum_{j=0}^3 \int_\Gamma K_{j, L,R}(\bx,\by) \sigma_j(\by) \,{\rm{d}}\by,
\end{equation}
where
\begin{equation}
\begin{aligned}
  K_{0,L,R}(\bx,\by) &= \partial_{y_1}^3 G_{L,R}(\bx,\by)  + 2 \partial_{y_1}\partial_{y_2}^2 G_{L,R}(\bx,\by) ,\\
  K_{1,L,R}(\bx,\by) &= \partial_{y_1}^2 G_{L,R}(\bx,\by)  + 2 \partial_{y_2}^2 G_{L,R}(\bx,\by) ,\\
  K_{2,L,R}(\bx,\by) &= \partial_{y_1} G_{L,R}(\bx,\by) ,\\
  K_{3,L,R}(\bx,\by) &= G_{L,R}(\bx,\by).
\end{aligned}
\end{equation}
Since our ansatz for $u_{L,R}$ is constructed from the appropriate domain fundamental solutions, they automatically satisfy the PDE in each half-space, the boundary conditions, and the radiations conditions away from $\Gamma$ (see \cite{epstein2023solvinga,epstein2023solvingb} for a thorough analysis of the Helmholtz case) for any $\{\sigma_j\}$.
Since $w_{L,R}$ are smooth, the jump relations for free-space kernels (see Appendix A of \cite{nekrasov2025boundary}) tell us that $u_{L,R}$ satisfy \eqref{eq:cont} provided $\{\sigma_j\}$ satisfy the integral equations
\begin{equation}\label{eq:matchIE}
\begin{pmatrix}
    -\sigma_0\\
    \sigma_1\\
    -\sigma_2\\
    \sigma_3
\end{pmatrix} + \begin{pmatrix}
    \cK_{0,0} & \cK_{0,1} & \cK_{0,2}&\cK_{0,3}\\
    \cK_{1,0} & \cK_{1,1} & \cK_{1,2}&\cK_{1,3}\\
    \cK_{2,0} & \cK_{2,1} & \cK_{2,2}&\cK_{2,3}\\
    \cK_{3,0} & \cK_{3,1} & \cK_{3,2}&\cK_{3,3}
\end{pmatrix}
\begin{pmatrix}
    \sigma_0\\
    \sigma_1\\
    \sigma_2\\
    \sigma_3
\end{pmatrix}=\begin{pmatrix}
    r_0\\
    r_1\\
    r_2\\
    r_3
\end{pmatrix}\quad \text{on } \Gamma,
\end{equation}
where $\cK_{i,j}$ are the integral operator whose kernels are
\begin{equation}
    K_{i,j}(\bx,\by) = \partial_{x_1}^iK_{j,R}(\bx,\by) - \partial_{x_1}^iK_{j,L}(\bx,\by).
\end{equation}
Since the kernels are all based on differences of domain Green's functions, the $G$ terms cancel. The kernels are therefore all smooth along $\Gamma$. Unfortunately, $w_{L,R}(\cdot,\by)$ decay algebraically like $|x_2|^{-1/2}e^{ik |x_2|}$ as $\bx\to\infty$ along $\Gamma,$ and the other derivatives decay similarly. It is therefore not practical to truncate the curve $\Gamma$ at any reasonable distance. To overcome this slow decay, we use the complex scaling method \cite{epstein2025complex,goodwill2024numerical,epstein2024coordinate,agocs2025complex}, which is based on the observation that $w_{L,R}$ are analytic and that if $x_2 = a+ (1+si)t$ for $a\in\bbR$ and $0<s<\infty$ then all of the kernels decay exponentially as $t\to\infty$ goes to infinity. Provided the right hand side is analytic, we may therefore analytically continue \eqref{eq:matchIE} to a contour $\tilde \Gamma$ parameterized by
\begin{equation}
    x_2(t) = t + 1i\lp \phi(t) - \phi(-t)\rp,
\end{equation}
where
\begin{equation}
    \phi(t)=
    \frac12(t+L)\operatorname{erfc}\lp3(t+L)\rp - \exp(-9(t+L)^2)/\sqrt{\pi}/6,
\end{equation}
which has slope $s=1$ at infinity. With this choice, $x_2(t)$ is real to machine precision in $[-L+2,L-2]$ and $L$ is chosen so that this interval contains all $x_2$ in the computational domain in \figref{fig:free_glue}. If the integrals defining $\cK_{i,j}$ are taken along $\tilde\Gamma$, then the integral equation will be Fredholm second kind on the set of continuous functions decaying exponentially along $\tilde \Gamma$ (see \cite{epstein2025complex}). If the data is analytic, then the $\sigma_j$-s will be analytic and agree with the densities that would be obtained by solving \eqref{eq:matchIE} on the real $\Gamma$ and so are the correct solution. Most physically meaningful data, such as data generated by a point source or trapped mode in either side will be analytic, and so the complex scaling method will give the correct solution. Once, the analytically continued integral equation has been solved, $u_{L,R}(\bx)$ can be computed using \eqref{eq:matchrep} with the integral taken along $\tilde\Gamma$. This value will be correct as long the contour deformation has not crossed the branch cut of $G_{L,R}$, which is inherited from $G$. This only occurs for $x_2$ outside of the region $\tilde \Gamma$ is real (i.e. when $|x_2|>L-2$). 

Since the kernels and densities decay exponentially along the deformed contour $\tilde\Gamma$, truncation gives errors which are exponentially small in the truncation height (see \cite{epstein2025complex}). To evaluate the layer potentials, we split $u_{L,R}$ into the contribution for the free-space Green's function, $u_{0,L,R}$, and the contribution of $w_{L,R}$, denoted $u_{w,L,R}$. Since $G$ is cheap to evaluate, $u_{0,L,R}$ can be computed quickly. We use adaptive integration to handle the near singularity. To evaluate the contribution of $w_{L,R}$, we let $K_{\xi,j,L,R}$ be the kernel that is equivalent to $K_{j,L,R}$ built out of $G_\xi$. The Floquet--Bloch transform of $u_{w,L,R}$ will be a solution of \eqref{eq:quaspde} with boundary data
\begin{equation}\label{eq:quasi_lump_BC}
\begin{aligned}
    f_{\xi,L,R} &= \mathcal{B}_{1,L,R}\left[\int_{\tilde\Gamma}\sum_j K_{\xi,j,L,R}(\bx,\by)\sigma_j(\by)\,d\by\right] \quad\text{and}\\ g_{\xi,L,R} &= \mathcal{B}_{2,L,R}\left[\int_{\tilde\Gamma}\sum_j K_{\xi,j,L,R}(\bx,\by)\sigma_j(\by)\,d\by\right],
\end{aligned}
\end{equation}
where $\mathcal{B}_{i,L,R}$ are the operators defining the boundary conditions on $\gamma_{L,R}$. Since $\gamma_{L,R}$ are separated from $\Gamma$, these can easily be evaluated using the smooth quadrature rule. The method described above can thus be used to efficiently evaluate $u_{w,L,R}$. The same method can be used to efficiently apply the system matrix of \eqref{eq:matchIE}, allowing us to efficiently solve \eqref{eq:matchIE} using GMRES \cite{GMRES}. Further details on the discretization of this are given in \cite{agocs2025complex}. Since each iteration requires an application of the dense matrices that appear in \eqref{eq:quasi_lump_BC} and the equivalent formula for evaluating the layer potentials supported on $\gamma_{L,R}$ at $\Gamma,$ we first precompute and store a low rank factorization of each of these matrices. This allows us to greatly accelerate each GMRES iteration.

\begin{remark}
    A nearly identical approach could be used to simulate scattering from the intersection of semi-infinite periodic gratings, as was considered in \cite{agocs2025complex}. The only difference is that there will be a singular interaction between the interface $\Gamma$ and the walls $\gamma_{L,R}$. To remove this singularity the representation~\eqref{eq:matchrep} would have to be modified to use the half-space Green's function (see Lemma 5 of~\cite{agocs2025complex}) and care would have to be taken to avoid catastrophic cancellations.
\end{remark}

\subsection{Numerical examples}

\begin{figure}
    \centering
    \includegraphics[width=0.6\linewidth]{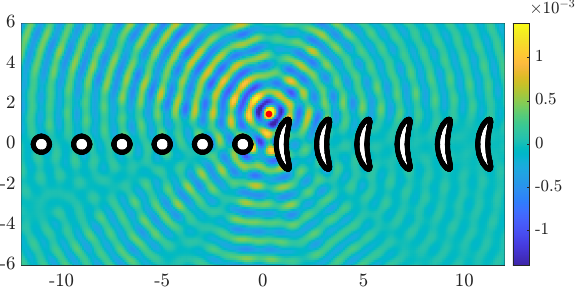}
    \caption{Imaginary part of the solution generated by a point source at (0.3,1.5) when free plate conditions are imposed on the black boundaries.}
    \label{fig:free_glue}
\end{figure}

As an illustration of the method, we consider free plate boundary conditions imposed on the curves $\gamma_{L,R}$ shown in \figref{fig:free_glue} and use the method described in the previous section to find the field generated by a point source at $\bx_0= (0.3,1.5) \in \Omega_R$. To do this, we let
\begin{equation}
    v(\bx) = \begin{cases}
        G_L(\bx,\bx_0) & \bx\in \Omega_L\\
        G_R(\bx,\bx_0)  & \bx\in \Omega_R
    \end{cases}.
\end{equation}
We then solve the equation \eqref{eq:matchIE} with $r_i = \partial_{x_1}^i G_L(\bx,\bx_0)-\partial_{x_1}^i G_R(\bx,\bx_0)$, which guarantees that
\begin{equation}
    u(\bx) = v(\bx)+\begin{cases}
        u_L(\bx) & \bx\in \Omega_L\\
        u_R(\bx) & \bx\in \Omega_R
    \end{cases}
\end{equation}
satisfies $\Delta^2 u -k^4 u = \delta(\bx-\bx_0)$ in $\Omega_{L,R}$ and has no jump in its first 3 derivatives across $\Gamma$. The resulting $u$ is shown in \figref{fig:free_glue}, which was computed in 18 GMRES iterations to converge when the tolerance was set to $10^{-10}$. 

\begin{figure}
    \centering
    \includegraphics[width=0.6\linewidth]{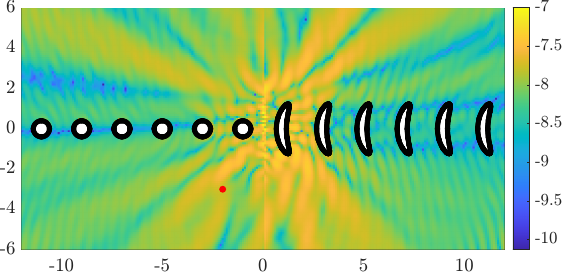}
    \caption{The error in our analytic solution test for the problem shown in \figref{fig:free_glue}. The color indicates the logarithm of the error divided by the maximum of $|v(\bx)|$ on the computational domain.}
    \label{fig:free_glue_error}
\end{figure}

In order to test the accuracy of the method, we repeat this experiment with
\begin{equation}
    v(\bx) = \begin{cases}
        0 & \bx\in \Omega_L\\
        G_R(\bx,\bx_0) & \bx\in \Omega_R
    \end{cases}.
\end{equation}
and $\bx_0=(-2,-3)\in \Omega_L$, and the $r_i$'s modified accordingly.
The resulting total field $u$ now satisfies $\Delta^2 u -k^4 u =0$ with zero boundary data, and hence is identically zero. The resulting error is shown in \figref{fig:free_glue_error}. As we can see, the solution is accurate to 7 digits everywhere. This computation took 20 GMRES iterations.

\section{Conclusion}\label{sec:conclusion}

In this work, we have described a numerical method for simulating time-harmonic flexural wave scattering from domains with infinite or semi-infinite periodic boundaries or obstacles. The approach is based on using a Floquet--Bloch transform to convert the problem to a set of quasi-periodic flexural scattering problems, which are then further reduced to integral equations on the boundary of the domain. The integral equation at each value of the Floquet--Bloch parameter $\xi$ is Fredholm second kind on the unit cell and can be efficiently solved by noting that the quasi-periodic Green's functions are equal to the free-space Green's function plus a smooth correction. Further, we showed how our method can be used to look for trapped modes, and to simulate scattering from junctions between two semi-infinite lines of obstacles. The method can be extended to handle junctions of semi-infinite gratings, with either clamped, supported, or free plate boundary conditions. Finally, in many applications the plate is floating on a fluid, e.g.\ an ice sheet floating on water. In principle, the approach considered here could be extended to these models, using as starting point the numerical methods given in\cite{surflay,flexgravwav}.

\section*{Conflict of Interest}
The authors declared that they have no conflict of interest.

\section*{Funding statement}
J.G.H. was supported in part by a Sloan Research Fellowship. This research was supported in part by grants from the NSF (DMS-2235451) and Simons Foundation (MPS-NITMB-00005320) to the NSF-Simons National Institute for Theory and Mathematics in Biology (NITMB).

\section*{Author Contributions}
All authors contributed to the conception of the paper. Code implementation and solution visualization were performed by Peter Nekrasov and Tristan Goodwill. The manuscript was written and edited by all of the authors.

\section*{Acknowledgements}
The authors would like to thank Douglas MacAyeal and Manas Rachh for many useful discussions.

\bibliography{references.bib}
\end{document}